\documentclass[reqno]{amsart}
\usepackage{amssymb,setspace}
\usepackage{ifpdf}
\ifpdf
 \usepackage[hyperindex,pagebackref]{hyperref}
\else
 \expandafter\ifx\csname dvipdfm\endcsname\relax
 \usepackage[hypertex,hyperindex,pagebackref]{hyperref}
 \else
 \usepackage[dvipdfm,hyperindex,pagebackref]{hyperref}
 \fi
\fi
\allowdisplaybreaks[4]
\numberwithin{equation}{section}
\theoremstyle{plain}
\newtheorem{thm}{Theorem}[section]

\newtheorem{lem}{Lemma}[section]
\theoremstyle{remark}
\newtheorem{rem}{Remark}[section]

\begin{document}

\title[Limits for ratios of polygamma functions at singularities]
{Limit formulas for ratios of polygamma functions at their singularities}

\author[F. Qi]{Feng Qi}
\address{Department of Mathematics, School of Science, Tianjin Polytechnic University, Tianjin City, 300387, China; School of Mathematics and Informatics, Henan Polytechnic University, Jiaozuo City, Henan Province, 454010, China}
\email{\href{mailto: F. Qi <qifeng618@gmail.com>}{qifeng618@gmail.com}, \href{mailto: F. Qi <qifeng618@hotmail.com>}{qifeng618@hotmail.com}, \href{mailto: F. Qi <qifeng618@qq.com>}{qifeng618@qq.com}}
\urladdr{\url{http://qifeng618.wordpress.com}}



\begin{abstract}
In the paper the author presents the limit formulas for ratios of polygamma functions at their singularities.
\end{abstract}

\thanks{This paper was typeset using \AmS-\LaTeX}

\maketitle

\section{Introduction}

Throughout this paper, we use $\mathbb{N}$ to denote the set of all positive integers.
\par
It is well known that the gamma function $\Gamma(z)$ is single valued and analytic over the entire complex plane, save for the points $z=-n$, with $n\in\{0\}\cup\mathbb{N}$, where it possesses simple poles with residue $\frac{(-1)^n}{n!}$. Its reciprocal $\frac1{\Gamma(z)}$ is an entire function possessing simple zeros at the points $z=-n$, with $n\in\{0\}\cup\mathbb{N}$. See related texts in~\cite[p.~255, 6.1.3]{abram}.
\par
It is also well known that the polygamma functions are defined by $\psi(z)=\frac{\Gamma'(z)}{\Gamma(z)}$ and $\psi^{(i)}(z)$ for $i\in\mathbb{N}$. Among them, the first five functions $\psi(z)$, $\psi'(z)$, $\psi''(z)$, $\psi^{(3)}(z)$, and $\psi^{(4)}(z)$ are known as the di-, tri-, tetra-, pentra-, and hexa-gamma functions respectively.
The polygamma function $\psi^{(n)}(z)$ for $n\in\{0\}\cup\mathbb{N}$ is single valued and analytic over the entire complex plane, save at the points $z=-m$, with $m\in\{0\}\cup\mathbb{N}$, where it possesses poles of order $n+1$. See related texts in~\cite[p.~260, 6.4.1]{abram}.
\par
In~\cite{Prabhu-arxiv, Prabhu-Srivastava}, the limit formulas
\begin{equation}\label{gamma-limit-eq}
\lim_{z\to-k}\frac{\Gamma(nz)}{\Gamma(qz)}=(-1)^{(n-q)k}\frac{q}{n}\cdot\frac{(qk)!}{(nk)!}
\end{equation}
and
\begin{equation}\label{polygamma-limit-eq}
  \lim_{z\to-k}\frac{\psi(nz)}{\psi(qz)}=\frac{q}{n}
\end{equation}
for any non-negative integer $k$ and all positive integers $n$ and $q$ were established.
\par
The main result of this paper may be stated as a theorem below.

\begin{thm}\label{polygamma-singularity-thm}
For $i,n,q\in\mathbb{N}$ and $k\in\{0\}\cup\mathbb{N}$, we have
\begin{equation}\label{deriv-polygamma-singul-lim}
\lim_{z\to-k}\frac{\psi^{(i)}(nz)}{\psi^{(i)}(qz)}= \biggl(\dfrac{q}{n}\biggr)^{i+1}.
\end{equation}
\end{thm}

\section{A lemma}

In~\cite[Theorem~1.2]{derivative-tan-cot.tex} the explicit formulas for the $n$-th derivatives of the cotangent function $\cot x$ was discovered. For proving the limit~\eqref{deriv-polygamma-singul-lim}, we cite and reorganize these formulas as Lemma~\ref{cotan-deriv-thm} below.

\begin{lem}\label{cotan-deriv-thm}
For $n\in\mathbb{N}$, the derivatives of the cotangent function may be computed by
\begin{equation}\label{cotan-deriv-2n-1}
\cot^{(2n-1)}x=\frac1{\sin^{2n}x}\sum_{i=0}^{n-1}b_{2n-1,2i}\cos(2ix)
\end{equation}
and
\begin{equation}\label{cotan-deriv-2n}
\cot^{(2n)}x=\frac1{\sin^{2n+1}x}\sum_{i=0}^{n-1}b_{2n,2i+1}\cos[(2i+1)x],
\end{equation}
where
\begin{equation}\label{a-1-0=1}
b_{1,0}=-1,
\end{equation}
\begin{equation}\label{2n+1=0=eq-thm}
b_{2n-1,0}=2n\sum_{\ell=0}^{n-2}(-1)^{\ell+1}\binom{2n-1}{\ell}(n-\ell-1)^{2n-2}
\end{equation}
for $n>1$,
\begin{equation}\label{a=2n+1=2i=thm}
b_{2n-1,2i}=2\sum_{\ell=0}^{n-i-1}(-1)^{\ell+1}\binom{2n}{\ell}(n-i-\ell)^{2n-1}
\end{equation}
for $1\le i\le n-1$, and
\begin{equation}\label{a=2n=2n-2i-1-use-thm}
b_{2n,2i+1}=2\sum_{\ell=0}^{n-i-1}(-1)^{\ell}\binom{2n+1}{\ell}(n-i-\ell)^{2n}
\end{equation}
for $0\le i\le n-1$.
\end{lem}

\begin{rem}
The equalities~\eqref{a=2n+1=2i=thm} and~\eqref{a=2n=2n-2i-1-use-thm} can be unified as
\begin{equation}
b_{p,q}=(-1)^{\frac{1-(-1)^p}2}2\sum_{\ell=0}^{\frac{p-q-1}2} (-1)^{\ell}\binom{p+1}{\ell}\biggl(\frac{p-q-1}2-\ell+1\biggr)^p
\end{equation}
for $0<q<p$.
\end{rem}

\section{Proof of Theorem~\ref{polygamma-singularity-thm}}

Now we set off to prove Theorem~\ref{polygamma-singularity-thm}.
\par
In~\cite[p.~260, 6.4.7]{abram}, the reflection formula
\begin{equation}\label{abram-6.4.7}
\psi^{(n)}(1-z)+(-1)^{n+1}\psi^{(n)}(z)=(-1)^n\cot^{(n)}(\pi z)
\end{equation}
is collected. Hence, we have
\begin{equation}\label{limit-general}
\lim_{z\to-k}\frac{\psi^{(i)}(nz)}{\psi^{(i)}(qz)}=\lim_{z\to-k}\frac{(-1)^i\pi\cot^{(i)}(\pi nz) -\psi^{(i)}(1-nz)}{(-1)^i\pi\cot^{(i)}(\pi qz) -\psi^{(i)}(1-qz)}.
\end{equation}
\par
When $i=1$, by~\eqref{limit-general} and Lemma~\ref{cotan-deriv-thm} applied to $n=1$, we have
\begin{align*}
\lim_{z\to-k}\frac{\psi'(nz)}{\psi'(qz)}&=\lim_{z\to-k}\frac{-\pi\cot'(\pi nz) -\psi'(1-nz)}{-\pi\cot'(\pi qz) -\psi'(1-qz)}\\
&=\lim_{z\to-k}\frac{\pi\cot'(\pi nz)+\psi'(1-nz)} {\pi\cot'(\pi qz) +\psi'(1-qz)}\\
&=\lim_{z\to-k}\frac{-\frac{\pi}{\sin^2(n\pi z)}+\psi'(1-nz)} {-\frac{\pi}{\sin^2(q\pi z)} +\psi'(1-qz)}\\
&=\lim_{z\to-k}\biggl[\frac{-{\pi}+{\sin^2(n\pi z)}\psi'(1-nz)} {-{\pi} +{\sin^2(q\pi z)}\psi'(1-qz)}\cdot\frac{\sin^2(q\pi z)}{\sin^2(n\pi z)}\biggr]\\
&=\lim_{z\to-k}\frac{\sin^2(q\pi z)}{\sin^2(n\pi z)}\\
&=\biggl(\frac{q}{n}\biggr)^2.
\end{align*}
When $i=2j$ and $j\in\mathbb{N}$, we have
\begin{multline*}
\begin{aligned}
\lim_{z\to-k}\frac{\psi^{(2j)}(nz)}{\psi^{(2j)}(qz)}&=\lim_{z\to-k}\frac{(-1)^{2j}\pi\cot^{(2j)}(\pi nz) -\psi^{(2j)}(1-nz)}{(-1)^{2j}\pi\cot^{(2j)}(\pi qz) -\psi^{(2j)}(1-qz)}\\
&=\lim_{z\to-k}\frac{\frac{\pi}{\sin^{2j+1}(nz\pi)}\sum_{i=0}^{j-1}b_{2j,2i+1}\cos[(2i+1)nz\pi]
 -\psi^{(2j)}(1-nz)}{\frac{\pi}{\sin^{2j+1}(qz\pi)}\sum_{i=0}^{j-1}b_{2j,2i+1}\cos[(2i+1)qz\pi] -\psi^{(2j)}(1-qz)}
\end{aligned}\\
\begin{aligned}
&=\lim_{z\to-k}\Biggl[\frac{{\pi}\sum_{i=0}^{j-1}b_{2j,2i+1}\cos[(2i+1)nz\pi]
 -{\sin^{2j+1}(nz\pi)}\psi^{(2j)}(1-nz)}{{\pi} \sum_{i=0}^{j-1}b_{2j,2i+1}\cos[(2i+1)qz\pi] -{\sin^{2j+1}(qz\pi)} \psi^{(2j)}(1-qz)}\cdot\frac{\sin^{2j+1}(qz\pi)} {\sin^{2j+1}(nz\pi)}\Biggr]\\
&=\lim_{z\to-k}\frac{\sin^{2j+1}(qz\pi)} {\sin^{2j+1}(nz\pi)}\\
&=\biggl(\frac{q}n\biggr)^{2j+1}.
\end{aligned}
\end{multline*}
When $i=2j+1$ and $j\in\mathbb{N}$, we have
\begin{multline*}
\begin{aligned}
\lim_{z\to-k}\frac{\psi^{(2j+1)}(nz)}{\psi^{(2j+1)}(qz)}&=\lim_{z\to-k}\frac{(-1)^{2j+1}\pi\cot^{(2j+1)}(\pi nz) -\psi^{(2j+1)}(1-nz)}{(-1)^{2j+1}\pi\cot^{(2j+1)}(\pi qz) -\psi^{(2j+1)}(1-qz)}\\
&=\lim_{z\to-k}\frac{\frac{\pi}{\sin^{2j+2}(nz\pi)}\sum_{i=0}^{j}b_{2j+1,2i}\cos(2inz\pi)
 +\psi^{(2j+1)}(1-nz)}{\frac{\pi}{\sin^{2j+2}(qz\pi)}\sum_{i=0}^{j}b_{2j+1,2i}\cos(2iqz\pi) +\psi^{(2j+1)}(1-qz)}
\end{aligned}\\
\begin{aligned}
&=\lim_{z\to-k}\Biggl[\frac{{\pi}\sum_{i=0}^{j}b_{2j+1,2i}\cos(2inz\pi)
 -{\sin^{2j+2}(nz\pi)}\psi^{(2j+1)}(1-nz)}{{\pi} \sum_{i=0}^{j}b_{2j+1,2i}\cos(2iqz\pi) -{\sin^{2j+2}(qz\pi)} \psi^{(2j+1)}(1-qz)}\cdot\frac{\sin^{2j+2}(qz\pi)} {\sin^{2j+2}(nz\pi)}\Biggr]\\
&=\lim_{z\to-k}\frac{\sin^{2j+2}(qz\pi)} {\sin^{2j+2}(nz\pi)}\\
&=\biggl(\frac{q}n\biggr)^{2j+2}.
\end{aligned}
\end{multline*}
In conclusion, the equality~\eqref{deriv-polygamma-singul-lim} for $i\ge1$ is proved. Combining~\eqref{polygamma-limit-eq} with the equality~\eqref{deriv-polygamma-singul-lim} for $i\ge1$ leads to Theorem~\ref{polygamma-singularity-thm}.

\end{document}